# Cohomological properties of the quantum shuffle product and application to the construction of quasi-Hopf algebras.

Cyrille Ospel

Laboratoire de Mathématiques, Université Blaise Pascal
Campus des Cézeaux, 63177 Aubière Cedex France
e-mail: ospel@ucfma.univ-bpclermont.fr

**Abstract:**

For a commutative algebra the shuffle product is a morphism of complexes. We generalize this result to the quantum shuffle product, associated to a class of non-commutative algebras (for example all the Hopf algebras). As a first application we show that the Hochschild-Serre identity is the dual statement of our result. In particular, we extend this identity to Hopf algebras. Secondly, we clarify the construction of a class of quasi-Hopf algebras.

## Introduction

The shuffles were defined by S. Eilenberg and S. MacLane to give an explicit formula for the equivalence of complexes of the Eilenberg-Zilber theorem. They were later used to show that the homology of an abelian group (or a commutative algebra) is an algebra with the shuffle product. On the other hand, the shuffles were used, more implicitly, by G. Hochschild and J.P. Serre [7] in the definition of the Hochschild-Serre identity of a group. Afterwards N. Habegger, V. Jones, O. Pino Ortiz, J. Ratcliffe [6] gave a formulation of the identity in terms of shuffles.

These two results were proved separately, by a long and technical verification of the two terms of the equality. In this paper, we prove that these results have a strong interaction. In fact, we show that the Hochschild-Serre identity is a consequence of the homological property of the shuffle map.

More generally, we consider algebras with an automorphism $\sigma$ of the square tensor product and some relations between the product and $\sigma$; such a couple is called a braided $\sigma$-commutative algebra. For such algebras, we show first that the quantum shuffle product, associated to $\sigma$ and defined by M. Rosso [16], can be factorized by the shuffle map. This factorization allows to prove that the quantum shuffle product is a morphism of complexes from a braided tensor product of chain complexes to the Hochschild chain complex of the algebra. The dual statement is a Hochschild-Serre identity for braided $\sigma$-commutative algebras.

A class of examples of such algebras is the class of Hopf algebras with invertible antipode. The braidind is given by the Woronowicz braiding [20]. In particular the case of the Hopf algebra of a group gives the classical Hochschild-Serre identity of the group. Secondly, we give a multiplicative statement of the Hochschild-Serre identity for a cocommutative Hopf algebra. Then, we use this result to clarify some constructions of quasi-Hopf algebras, associated to the Drinfeld Double, and defined by R. Dijkgraaf, V. Pasquier, P. Roche [3] and D. Bulacu, F. Panaite [2].

In [14] we show that the homological property of the quantum shuffle product allows to extend the first iteration of the abelian group homology construction [5] to non-commutative Hopf algebras. The multiplicative cohomology associated to this chain complex has applications to the theory of invariants for links and 3-manifolds. For example, the 3-cocycles of this cohomology are weight systems for links.

### Acknowledgments

I would like to express my gratitude to Professor M. Rosso for explaining me the subject, by helpful discussions and bibliographic references, and for suggesting me the leading idea.

**Notations:** $\mathbb{K}$ is a commutative field.

All algebras considered are associative algebras over $\mathbb{K}$ with unit. The product is denoted by $\mu$.

We use Sweedler's notation for coproduct $\Delta(a) = \sum a_{(1)} \otimes a_{(2)}$.

Let $V$ be a vector space over $\mathbb{K}$. The tensor vector space $T(V)$ of $V$ is defined by $T(V) = \bigoplus_{k \geq 0} V^{\otimes k}$.

Let $\underline{w} \in V^{\otimes n}$. We denote the degree $n$ of $\underline{w}$ by $|\underline{w}|$.

$\Sigma_n$ is the set of all permutations of $\{1, \ldots, n\}$. For all $\nu$ in $\Sigma_n$, we denote the sign of $\nu$ by $(-1)^{|\nu|}$.

## 1 A new construction of the quantum shuffle product.

The original definition of the quantum shuffle product was given in the framework of representation theory. This point of view is not useful to study the homological properties of the product. So we will give a factorization of this product by morphisms of complexes, in particular by the shuffle map.

### 1.1 The original definition.

The quantum shuffle product was first defined by M. Rosso [15]. It describes the product of the following cotensorial Hopf algebra. Let $H$ be a Hopf algebra and $M$ an $H$-Hopf bimodule, with the bicomodule structure given by $\delta_L$ and $\delta_R$. The cotensorial

Hopf algebra $T_H^c(M)$ was defined by W. Nichols [13], by: $T_H^c(M) = H \oplus \bigoplus_{n \geq 1} M^{\square n}$, where $M \square M$ is the kernel of $\delta_R \otimes \text{Id}_M - \text{Id}_M \otimes \delta_L$.

This Hopf algebra is an $H$-Hopf bimodule. As algebra $T_H^c(M)$ is the crossed product of $H$ by the left-coinvariant subspace of the cotensorial Hopf bimodule. The algebra structure induced on the left-coinvariant subspace is given by the quantum shuffle product.

More explicitly, let $V$ be a vector space over $\mathbb{K}$ and $\sigma \in \text{End}(V \otimes V)$ which satisfies the braid equation:

$$\sigma_2 \sigma_1 \sigma_2 = \sigma_1 \sigma_2 \sigma_1, \tag{1}$$

where for $n, i$ non negative integers such that $n > i$, $\sigma_i \in \text{End}(V^{\otimes n})$ is defined by $\sigma_i = \text{Id}^{\otimes i-1} \otimes \sigma \otimes \text{Id}^{\otimes n-i-1}$.

For $n \geq 2$, we denote by $T_\sigma$ the representation of the braid group $B_n$ on $V^{\otimes n}$ defined on the generators $(\omega_i)_{1 \leq i \leq n-1}$ by:

$$T_\sigma(\omega_i) = \sigma_i, \quad 1 \leq i \leq n-1.$$

For all non negative integers $p, q$ such that $p + q = n$, $S_{p,n}$ is the set of all $(p,q)$-shuffles, i.e. the set of all $w \in \sum_n$ such that $w(1) < \cdots < w(p)$ and $w(p+1) < \cdots < w(n)$. Let $(\nu_i)_{1 \leq i \leq n-1}$ be the transpositions $(i, i+1)$ of $\sum_n$. For $w \in \sum_n$ and $w = \nu_{i_1} \ldots \nu_{i_r}$ a reduced decomposition of $w$ we define the extension of $T_\sigma$ to $\Sigma_n$ by:

$$T_\sigma(w) = \sigma_{i_1} \ldots \sigma_{i_r}.$$

The quantum shuffle product $\varphi_\sigma : T(V) \otimes T(V) \to T(V)$, associated to $\sigma$, is defined for all positive integers $p + q = n$ by:

$$\varphi_\sigma : \begin{array}{rcl} V^{\otimes p} \otimes V^{\otimes q} & \longrightarrow & V^{\otimes n} \\ \underline{v} \otimes \underline{v}' & \longmapsto & \displaystyle\sum_{w \in S_{p,n}} T_\sigma(w)(\underline{v} \otimes \underline{v}'). \end{array}$$

If $p = 0$ or $q = 0$, the product is just the multiplication by elements of $\mathbb{K}$.

With this product and the unit $1_\mathbb{K}$, $T(V)$ is an algebra.

In this paper, we will always use a quantum shuffle product with sign, i.e. associated to the braid $-\sigma$. We will still denoted it by $\varphi_\sigma$. It is defined for all positive integers $p + q = n$ by:

$$\varphi_\sigma : \begin{array}{rcl} V^{\otimes p} \otimes V^{\otimes q} & \longrightarrow & V^{\otimes n} \\ \underline{v} \otimes \underline{v}' & \longmapsto & \displaystyle\sum_{w \in S_{p,n}} (-1)^{|w|} T_\sigma(w)(\underline{v} \otimes \underline{v}'). \end{array}$$

### 1.2 A twisted version of the shuffle map.

The shuffle map was defined by S. Eilenberg and S. MacLane [5] as a map from the tensor product of two simplicial chain complexes to the cartesian product of these chain complexes (for the classical definitions and results of the simplicial theory, we refer to S. MacLane [12] and to J.L. Loday [11]).

They have shown that this map is a morphism of complexes. The proof is based only on the properties of the simplicial objects and of the $(p,q)$-shuffles.

We use this map on the un-normalized bar resolution of an algebra $A$, defined by S. MacLane [12] as the following simplicial object $\beta(A,A)$:

$$\beta_n(A,A) = A \otimes A^{\otimes n} \otimes A,$$

and for $\underline{a} = a^0[a^1 \mid \cdots \mid a^n]a^{n+1} \in \beta_n(A,A)$ the operators are given by:

$$\begin{aligned}
d_0(\underline{a}) &= a^0 a^1[a^2 \mid \cdots \mid a^n]a^{n+1}, \\
d_i(\underline{a}) &= a^0[a^1 \mid \cdots \mid a^i a^{i+1} \mid \cdots \mid a^n]a^{n+1} & \forall\, i \in [1,..,n-1], \\
d_n(\underline{a}) &= a^0[a^1 \mid \cdots \mid a^{n-1}]a^n a^{n+1}, \\
s_i(\underline{a}) &= a^0[a^1 \mid \cdots \mid a^i \mid 1_A \mid a^{i+1} \mid \cdots \mid a^n]a^{n+1} & \forall\, i \in [0,..,n].
\end{aligned}$$

For the standard tensor product and cartesian product of simplicial chain complexes, S. Eilenberg and S. MacLane have proved:

**Theorem 1.1** *The shuffle map* $g : \beta(A,A) \otimes \beta(A,A) \to \beta(A,A) \times \beta(A,A)$ *,defined for all non negative integers* $p+q=n$ *and* $\underline{a} \in \beta_p(A,A)$, $\underline{b} \in \beta_q(A,A)$ *by*

$$g_n(\underline{a} \otimes \underline{b}) = \sum_{\nu \in S_{p,n}} (-1)^{|\nu|} s_{\nu(p+q)-1} \circ \cdots \circ s_{\nu(p+1)-1}(\underline{a}) \otimes s_{\nu(p)-1} \circ \cdots \circ s_{\nu(1)-1}(\underline{b}),$$

*is a morphism of complexes.*

We will twist this map by a braiding $\sigma$, i.e. an automorphism of $A \otimes A$. For this, we define a map:

$$\psi : \beta(A,A) \times \beta(A,A) \to \beta(A \otimes A, A \otimes A),$$

of degree 0. For $n \in \mathbb{N}$, the map $\psi_n$ is defined on $\beta_n(A,A) \otimes \beta_n(A,A)$ by:

$$\psi_n = \sigma_{2n+2} \circ (\sigma_{2n}\sigma_{2n+1}) \circ \cdots \circ (\sigma_4 \ldots \sigma_{n+3}) \circ (\sigma_2 \ldots \sigma_{n+2}).$$

In the classical case $\sigma$ is the flip $\tau$ defined on $A \otimes A$ by $\sigma(a \otimes b) = b \otimes a$ and the map $\psi$ is a morphism of complexes.

More generally, it is easy to see that for the classical structure of algebra on $A \otimes A$ and $\sigma$ different from $\tau$, $\psi$ is not a morphism of complexes.

Thus we braid the product of $A \otimes A$ by $\sigma$. Therefore, using the ideas of J. Baez [1], A. Joyal and R. Street [9], A. Van Daele and S. Van Keer [18], M. Wambst [19] we give conditions on $\sigma$ to define a braided algebra structure on $A \otimes A$.

**Proposition 1.1** *Let* $\mu_\sigma : (A \otimes A) \otimes (A \otimes A) \to A \otimes A$ *be the linear map defined by:*

$$\mu_\sigma = (\mu \otimes \mu)\sigma_2^{-1}.$$

*The vector space $A \otimes A$ is an algebra with product $\mu_\sigma$ and unit $1_\sigma = 1_A \otimes 1_A$ if and only if*

$$\sigma_1(\mu \otimes \mathrm{Id}) = (\mathrm{Id} \otimes \mu)\sigma_1\sigma_2, \quad \sigma_1(\mathrm{Id} \otimes \mu) = (\mu \otimes \mathrm{Id})\sigma_2\sigma_1 \qquad (2)$$

$$\sigma(a \otimes 1_A) = 1_A \otimes a, \quad \sigma(1_A \otimes a) = a \otimes 1_A \qquad \forall a \in A \qquad (3)$$

*This algebra is denoted by $A \otimes_\sigma A$.*

**Example:** For $\sigma$ the flip $\tau$ defined on $a, b \in A$ by $\tau(a \otimes b) = b \otimes a$, the braided algebra $A \otimes_\tau A$ is just the classical algebra $A \otimes A$.

With this braided algebra structure, we have:

**Proposition 1.2** *The map $\psi$ is a morphism of complexes between $\beta(A, A) \times \beta(A, A)$ and $\beta(A \otimes_\sigma A, A \otimes_\sigma A)$.*

**Proof:** Using the definition of the un-normalized bar resolution the proposition is proved if, for all integers $n, i$ such that $n > 0$ and $0 \leq i \leq n$, we show that:

$$\psi_{n-1} \circ (d_i^n \otimes d_i^n) = (\mathrm{Id}^{\otimes 2i} \otimes \mu_\sigma \otimes \mathrm{Id}^{\otimes 2n-2i}) \circ \psi_n.$$

To make this proof more intuitive, we use the graphical technique introduced by many authors, for example see C. Kassel [10].

The product $\mu$ of $A$ is represented by 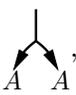, the map $\sigma$ by 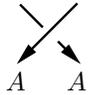, the map $\sigma^{-1}$ by 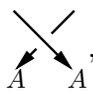, and the composition $f \circ g$ of two maps by 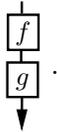.

So the map $\psi_n$ is represented by:

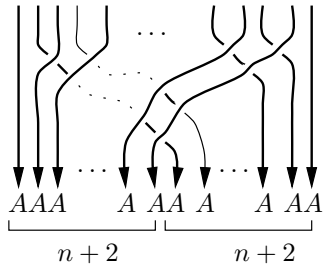

Now using the property (2), we prove the equality. For example, for $i = n$ the graphical proof is:

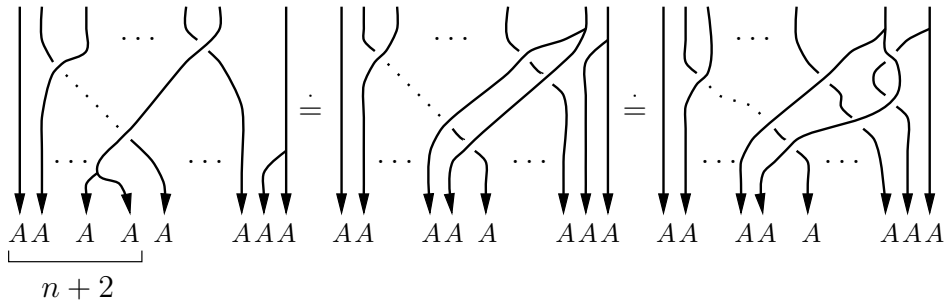

□

Thus, the map $\psi \circ g$ is a morphism of complexes from $\beta(A,A) \otimes \beta(A,A)$ to $\beta(A\otimes_\sigma A, A\otimes_\sigma A)$.

## 1.3 Braided bimodule structure and braided differential.

Let $M$ be an $A\otimes_\sigma A$-bimodule; we show that the map $\psi \circ g$ can be extended to a map of complexes from $M \otimes T(A) \otimes T(A)$ to $M \otimes T(A\otimes_\sigma A)$. For this, we first define a braided $A\otimes_\sigma A$-bimodule structure on $\beta(A,A) \otimes \beta(A,A)$, and then we produce a braided differential on $M \otimes T(A) \otimes T(A)$.

On the one hand, for an algebra $A$, the vector space $\beta(A,A)$ has a natural structure of $A$-bimodule: the left action, respectively right action, is given by the product with the first element, respectively the last element.
In particular let us consider the algebra $A\otimes_\sigma A$. Then using the associativity of the algebra $A\otimes_\sigma A$, one can prove that the vector space $\beta(A\otimes_\sigma A, A\otimes_\sigma A)$ is a chain complex of $A\otimes_\sigma A$-bimodules, i.e. the differential $d_\sigma$ of $\beta(A\otimes_\sigma A, A\otimes_\sigma A)$ is a morphism of $A\otimes_\sigma A$-bimodules.
Denote by $\widetilde{\otimes}$ the tensor product of $A\otimes_\sigma A$-bimodules (where the tensor product of two $A$-bimodules $M$ and $N$ is defined by $M \otimes_{A\otimes A^{op}} N$).
The map $\mathrm{Id}_M \widetilde{\otimes} d_\sigma$ endows the vector space $M\widetilde{\otimes}\beta(A\otimes_\sigma A, A\otimes_\sigma A)$ with a structure of chain complex.
To extend on the left the map $\psi \circ g$, we use the injection:

$$\Gamma : M \otimes T(A \otimes_\sigma A) \to M\widetilde{\otimes}\beta(A\otimes_\sigma A, A\otimes_\sigma A),$$

given by $\Gamma(\alpha \otimes \underline{w}) = \alpha \otimes 1_\sigma[\underline{w}]1_\sigma$.
Now using the precedent bimodule structure, we have :

**Lemma 1.1** *With the classical structure of Hochschild complex on $M \otimes T(A\otimes_\sigma A)$, the injection $\Gamma$ is an isomorphism of complexes.*

**Remark:** This result is not specific to $A\otimes_\sigma A$. In fact, Mac Lane [12] has shown that for any algebra $A_1$ and any $A_1$-bimodule $M_1$ the corresponding injection from $M_1 \otimes T(A_1)$ to $M_1\widetilde{\otimes}\beta(A_1, A_1)$ is an isomorphism of complexes.

On the other hand, the natural $A$-bimodule structure of $\beta(A,A)$ can be extended to a $A\otimes_\sigma A$-bimodule structure on $\beta(A,A) \otimes \beta(A,A)$. For all non negative integers $p,q$, the left action $\varphi_L$ and the right action $\varphi_R$ are defined by:

$$\begin{aligned}\varphi_L &: A\otimes_\sigma A \otimes \beta_p(A,A) \otimes \beta_q(A,A) \longrightarrow \beta_p(A,A) \otimes \beta_q(A,A) \\ \varphi_L &= (\mu \otimes \mathrm{Id}^{\otimes p+1} \otimes \mu \otimes \mathrm{Id}^{\otimes q+1}) \circ \sigma^{-1}_{p+3}\ldots\sigma^{-1}_2,\end{aligned}$$

$$\begin{aligned}\varphi_R &: \beta_p(A,A) \otimes \beta_q(A,A) \otimes A \otimes_\sigma A \longrightarrow \beta_p(A,A) \otimes \beta_q(A,A) \\ \varphi_R &= (\mathrm{Id}^{\otimes p+1} \otimes \mu \otimes \mathrm{Id}^{\otimes q+1} \otimes \mu) \circ \sigma^{-1}_{p+3}\ldots\sigma^{-1}_{p+q+4}.\end{aligned}$$

The graphical representations of these actions are:

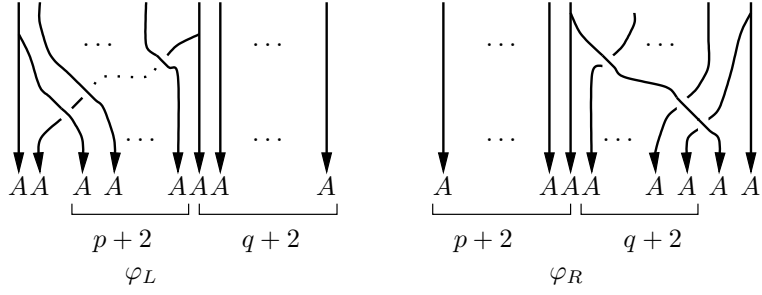

$$\varphi_L \qquad\qquad\qquad \varphi_R$$

Using the graphical technique one can show that with these structures the vector spaces $\beta(A,A) \otimes \beta(A,A)$ and $\beta(A,A) \times \beta(A,A)$ are $A \otimes_\sigma A$-bimodules and that the vector space $\beta(A,A) \otimes \beta(A,A)$ is a chain complex of $A \otimes_\sigma A$-bimodules. In particular, for all non negative integer $i$, the maps $\mathrm{Id} \otimes d_i$ and $d_i \otimes \mathrm{Id}$ are morphisms of $A \otimes_\sigma A$-bimodules.

So, for $d_\otimes$ the differential of $\beta(A,A) \otimes \beta(A,A)$, the map $\mathrm{Id}_M \widetilde{\otimes} d_\otimes$ endows $M \widetilde{\otimes} (\beta(A,A) \otimes \beta(A,A))$ with a structure of chain complex.

Now using these structures, we have:

**Lemma 1.2** *For all $p, q \in \mathbb{N}$, the injection:*

$$\Omega_{p,q} : M \otimes A^{\otimes p} \otimes A^{\otimes q} \to M \widetilde{\otimes} (\beta_p(A,A) \otimes \beta_q(A,A))$$

*given by $\Omega_{p,q}(\alpha \otimes \underline{w} \otimes \underline{w}') = \alpha \otimes 1_A[\underline{w}]1_A \otimes 1_A[\underline{w}']1_A$, induces an isomorphism of complexes. The differential $\rho$ on $M \otimes T(A) \otimes T(A)$ is defined, for $n \geq 1$, by $\rho_n = \displaystyle\sum_{p+q=n} \rho_{p,q}$, where:*

$$\rho_{p,q} : M \otimes A^{\otimes p} \otimes A^{\otimes q} \longrightarrow \left[M \otimes A^{\otimes p-1} \otimes A^{\otimes q}\right] \oplus \left[M \otimes A^{\otimes p} \otimes A^{\otimes q-1}\right],$$

*with, for $\underline{a} = a^1 \otimes \cdots \otimes a^p \in A^{\otimes p}$, $\underline{b} = b^1 \otimes \cdots \otimes b^q \in A^{\otimes q}$, and*

$$\begin{array}{rcll}
\sigma_q \ldots \sigma_1 : & A \otimes A^{\otimes q} & \to & A^{\otimes q} \otimes A \\
& a^p \otimes \underline{b} & \mapsto & \sum_i \underline{b}_i \otimes a_i^p,
\end{array} \qquad
\begin{array}{rcll}
\sigma_1 \ldots \sigma_p : & A^{\otimes p} \otimes A & \to & A \otimes A^{\otimes p} \\
& \underline{a} \otimes b^1 & \mapsto & \sum_j b_j^1 \otimes \underline{a}_j,
\end{array}$$

*we have*

$$\begin{aligned}
\rho_{p,q}(\alpha \otimes \underline{a} \otimes \underline{b}) = \ & \alpha.(a^1 \otimes 1_A) \otimes a^2 \otimes \cdots \otimes a^p \otimes \underline{b} \\
& + \sum_{i=1}^{p-1} (-1)^i \alpha \otimes a^1 \otimes \cdots \otimes a^i a^{i+1} \otimes \cdots \otimes a^p \otimes \underline{b} \\
& + (-1)^p \sum_i (a_i^p \otimes 1_A).\alpha \otimes a^1 \otimes \cdots \otimes a^{p-1} \otimes \underline{b}_i \\
& + (-1)^p \sum_j \alpha.(1_A \otimes b_j^1) \otimes \underline{b}_j \otimes b^2 \otimes \cdots \otimes b^q \\
& + \sum_{i=1}^{q-1} (-1)^{i+p} \alpha \otimes \underline{a} \otimes b^1 \otimes \cdots \otimes b^i b^{i+1} \otimes \cdots \otimes b^q \\
& + (-1)^{p+q} (1_A \otimes b^q).\alpha \otimes \underline{a} \otimes b^1 \otimes \cdots \otimes b^{q-1}.
\end{aligned}$$

**Proof:** Let $a'[\underline{a}]a'' \in \beta_p(A,A)$ and $b'[\underline{b}]b'' \in \beta_q(A,A)$. Then with the $A \otimes_\sigma A$-bimodule structure, we have:

$$a'[\underline{a}]a'' \otimes b'[\underline{b}]b'' = \sum_i (a' \otimes w_i) \cdot (1_A[\underline{w}_i]1_A \otimes 1_A[\underline{\tau}_i]1_A) \cdot (\tau_i \otimes b''),$$

where $\sigma_{p+q+1}\ldots\sigma_{p+2}\sigma_1\ldots\sigma_{p+1}: A^{\otimes p} \otimes A \otimes A \otimes A^{\otimes q} \to A \otimes A^{\otimes p} \otimes A^{\otimes q} \otimes A$
$$\underline{a} \otimes a'' \otimes b' \otimes \underline{b} \mapsto \sum_i w_i \otimes \underline{w}_i \otimes \underline{\tau}_i \otimes \tau_i.$$

So, the inverse of $\Omega_{p,q}$ is given by:

$$\Omega_{p,q}^{-1}(\alpha \otimes a'[\underline{a}]a'' \otimes b'[\underline{b}]b'') = \sum_i (\tau_i \otimes b'').\alpha.(a' \otimes w_i) \otimes \underline{w}_i \otimes \underline{\tau}_i.$$

Hence, using $\Omega$ we transport the differential $\operatorname{Id}_M \widetilde{\otimes} d_\otimes$ to the differential $\rho$. □

Therefore, to extend $\psi \circ g$ to a morphism of complexes from $M \otimes T(A) \otimes T(A)$ to $M \otimes T(A \otimes_\sigma A)$, we only have to prove that the maps $\operatorname{Id} \otimes g$ and $\operatorname{Id} \otimes \psi$ are well defined. For this, we must show that $g$ and $\psi$ are morphisms of $A \otimes_\sigma A$-bimodules.
This is just a straightforward verification (using the equality (3) and the graphical technique). In particular, for all non negative integer $i$, the maps $\operatorname{Id} \otimes s_i$ and $s_i \otimes \operatorname{Id}$ are morphisms of $A \otimes_\sigma A$-bimodules.

Thus, we have the theorem:

**Theorem 1.2** *Let $A$ be an algebra, $\sigma \in \operatorname{Aut}(A \otimes A)$ satisfying (2) and (3) and $M$ an $A \otimes_\sigma A$-bimodule.*
*There exists a map $\Theta$ from $T(A) \otimes T(A)$ to $T(A \otimes_\sigma A)$ such that $\operatorname{Id}_M \otimes \Theta$ is the morphism of complexes which extends $\psi \circ g$.*

**Proof:** We have shown that the map $\widetilde{\Theta} = \Gamma^{-1} \circ (\operatorname{Id}_M \widetilde{\otimes} \psi) \circ (\operatorname{Id}_M \widetilde{\otimes} g) \circ \Omega$ from $M \otimes T(A) \otimes T(A)$ to $M \otimes T(A \otimes_\sigma A)$ is a morphism of complexes.
So we only have to proof that for all non negative integer $n$, the map $\widetilde{\Theta}_n$ can be written as $\widetilde{\Theta}_n = \operatorname{Id}_M \otimes \Theta_n$, where:

$$\Theta_n : \bigoplus_{p+q=n} A^{\otimes p} \otimes A^{\otimes q} \to (A \otimes_\sigma A)^{\otimes n}.$$

For this, associate to any non negative integers $p, n, i$, such that $0 \leq p \leq n$, $0 \leq i \leq p$, the maps $s'_i : A^{\otimes p} \to A^{\otimes p+1}$ and $g' : A^{\otimes p} \otimes A^{\otimes n-p} \to A^{\otimes n} \otimes A^{\otimes n}$ by:

$$s'_i(a^1 \otimes \cdots \otimes a^p) = a^1 \otimes \cdots \otimes a^i \otimes 1_A \otimes a^{i+1} \otimes \cdots \otimes a^p,$$
$$g' = \sum_{\nu \in S_{p,n}} (-1)^{|\nu|} s'_{\nu(n)-1} \ldots s'_{\nu(p+1)-1} \otimes s'_{\nu(p)-1} \ldots s'_{\nu(1)-1}.$$

Then the map $\Theta$ is defined for all non negative integer $n$ by $\Theta_n = \psi_{n-2} \circ g'$. It is easy to show that $\widetilde{\Theta} = \operatorname{Id}_M \otimes \Theta$. □

### 1.4 Factorization of the quantum shuffle product.

We give a factorization of the quantum shuffle product by the map $\Theta$.

**Theorem 1.3** *Let $A$ be an algebra and let $\sigma \in \mathrm{Aut}(A \otimes A)$ satisfying (1), (2) and (3). The quantum shuffle product $\varphi_\sigma$ can be factorized by $\Theta$:*

$$\varphi_\sigma = \mu^\otimes \circ \Theta,$$

*where for all non negative integer $n$, $\mu^{\otimes n}$ is the product of $A$ tensorized $n$ times.*

The theorem is a consequence of the following proposition:

**Proposition 1.3** *For all non negative integers $p, n$ such that $n \geq 2$, $0 \leq p \leq n$ and $\nu \in S_{p,n}$, we have:*

$$T_\sigma(\nu) = \mu^{\otimes n} \circ \psi_{n-2} \circ (s'_{\nu(n)-1} \cdots s'_{\nu(p+1)-1} \otimes s'_{\nu(p)-1} \cdots s'_{\nu(1)-1}).$$

**Proof:** The proposition is clear if $p = 0$ or $p = n$. Thus assume that $p$ is such that $1 \leq p \leq n - 1$.

Then the proposition is proved by induction on $n$. For $n = 2$ it is clear.

Assume that the proposition is true for $n - 1$. Let $p, q \in \mathbb{N}$ be such that $p + q = n$, $1 \leq p \leq n - 1$ and $\nu \in S_{p,n}$. We have two possibilities: either $\nu(1) = 1$ or $\nu(p+1) = 1$. In the first case we define a new shuffle $\nu' \in S_{p-1, n-1}$ by:

$$\nu'(i) = \nu(i+1) - 1, \qquad 1 \leq i \leq n - 1.$$

For $a^1 \otimes \cdots \otimes a^p \otimes \underline{b} \in A^{\otimes p} \otimes A^{\otimes q}$ we have:

$$\begin{aligned}
&s'_{\nu(n)-1} \cdots s'_{\nu(p+1)-1}(a^1 \otimes \cdots \otimes a^p) \otimes s'_{\nu(p)-1} \cdots s'_{\nu(1)-1}(\underline{b}) \\
&= s'_{\nu(n)-1} \cdots s'_{\nu(p+1)-1}(a^1 \otimes \cdots \otimes a^p) \otimes s'_{\nu(p)-1} \cdots s'_{\nu(2)-1}(1_A \otimes \underline{b}) \\
&= a^1 \otimes s'_{\nu(n)-2} \cdots s'_{\nu(p+1)-2}(a^2 \otimes \cdots \otimes a^p) \otimes 1_A \otimes s'_{\nu(p)-2} \cdots s'_{\nu(2)-2}(\underline{b}) \\
&= a^1 \otimes s'_{\nu'(n-1)-1} \cdots s'_{\nu'(p)-1}(a^2 \otimes \cdots \otimes a^p) \otimes 1_A \otimes s'_{\nu'(p-1)-1} \cdots s'_{\nu'(1)-1}(\underline{b})
\end{aligned}$$

Using the induction and (3), we deduce the proposition in the first case.

In the second case, we define a shuffle $\nu' \in S_{p, n-1}$ by:

$$\begin{cases} \nu'(i) = \nu(i) - 1 & 1 \leq i \leq p, \\ \nu'(i) = \nu(i+1) - 1 & p+1 \leq i \leq n - 1. \end{cases}$$

For $\underline{a} \otimes b^1 \otimes \cdots \otimes b^q \in A^{\otimes p} \otimes A^{\otimes q}$ we have:

$$\begin{aligned}
&s'_{\nu(n)-1} \cdots s'_{\nu(p+1)-1}(\underline{a}) \otimes s'_{\nu(p)-1} \cdots s'_{\nu(1)-1}(b^1 \otimes \cdots \otimes b^q) \\
&= 1_A \otimes s'_{\nu(n)-2} \cdots s'_{\nu(p+2)-2}(\underline{a}) \otimes b^1 \otimes s'_{\nu(p)-2} \cdots s'_{\nu(1)-2}(b^2 \otimes \cdots \otimes b^q) \\
&= 1_A \otimes s'_{\nu'(n-1)-1} \cdots s'_{\nu'(p+1)-1}(\underline{a}) \otimes b^1 \otimes s'_{\nu'(p)-1} \cdots s'_{\nu'(1)-1}(b^2 \otimes \cdots \otimes b^q)
\end{aligned}$$

Using the induction and (3) we deduce the proposition in the second case. $\square$

**Remark:** If we define the quantum shuffle product by the formula of Theorem 1.3, we do not need to ask for $\sigma$ to satisfy (1). But this equality is necessary to prove the associativity.

## 2 Homological corollaries.

We give immediate corollaries of the factorization, which extend some classical results to braided algebras.

### 2.1 The quantum shuffle product as a morphism of complexes.

Let $M$ be an $A$-bimodule; if $\sigma$ satisfies

$$\mu \circ \sigma = \mu, \tag{4}$$

then $M$ is an $A \otimes_\sigma A$-bimodule with the structure defined, for all $a, b \in A$ and $w \in M$, by:

$$(a \otimes b).w = ab.w, \qquad w.(a \otimes b) = w.ab\ .$$

With this structure and Theorem (1.3), we deduce immediately the following corollary:

**Corollary 2.1** *Let $A$ be an algebra, $\sigma \in \mathrm{Aut}(A \otimes A)$ satisfying (2), (3), (4) and $M$ an $A$-bimodule.*
*The morphism $\mathrm{Id}_M \otimes \varphi_\sigma$ is a morphism of complexes, between the chain complex $(M \otimes T(A) \otimes T(A), \rho_\cdot)$ and the Hochschild chain complex of $A$ with coefficients in $M$.*

**Remark:** Such a pair $(A, \sigma)$ is called a braided $\sigma$-commutative algebra. All Hopf algebras with invertible antipode (in particular all finite dimensional Hopf algebras) are braided $\sigma$-commutative algebras, with $\sigma$ given by the Woronowicz braiding. This braiding, defined in [20], is given for $a, b \in A$ by $\sigma(a \otimes b) = \sum b_{(1)} \otimes \mathrm{ad}_{b_{(2)}} a$, with $\mathrm{ad}_a(b) = \sum \mathrm{S}(a_{(1)}) b a_{(2)}$ the adjoint action.
The inverse is given for $a, b \in A$ by $\sigma^{-1}(a \otimes b) = \sum a_{(3)} b \mathrm{S}^{-1}(a_{(2)}) \otimes a_{(1)}$.

We have the same result for the normalized chain complexes.
Let $I$ be the sub-space of $T(A)$ generated by the elements $a^1 \otimes \cdots \otimes a^n \in A^{\otimes n}$ where $n \in \mathbb{N}^*$ and one of the $a^i$ is equal to $1_A$.
It is well know that the Hochschild differential induces a differential on the quotients. Using property (3) of $\sigma$, it is easy to show that the differential $\rho$ and the product $\mathrm{Id}_M \otimes \varphi_\sigma$ induce respectively a differential and a product on $M \otimes T(A)/I \otimes T(A)/I$. So we have:

**Corollary 2.2** *Let $A$ be an algebra, $\sigma \in \mathrm{Aut}(A \otimes A)$ satisfying (2), (3), (4) and $M$ an $A$-bimodule.*
*The morphism $\mathrm{Id}_M \otimes \varphi_\sigma$ is a morphism of complexes, between the chain complex $(M \otimes T(A)/I \otimes T(A)/I, \rho_\cdot)$ and the normalized Hochschild chain complex of $A$ with coefficients in $M$.*

## 2.2 Braided square tensor product of a Hochschild complex with coefficients in $\mathbb{K}$.

Let $A$ be a braided $\sigma$-commutative algebra augmented with $\varepsilon : A \to \mathbb{K}$. Then $\mathbb{K}$ is an $A$-bimodule via the augmentation. Using Corollary (2.1), the product $\varphi_\sigma$ is a morphism of complexes between $(T(A) \otimes T(A), \rho)$ and the Hochschild complex of $A$ with coefficients in $\mathbb{K}$. The differential $\rho$ is defined for all non negative integers $p, q$ and $\underline{a} = a^1 \otimes \cdots \otimes a^p \in A^{\otimes p}$, $\underline{b} = b^1 \otimes \cdots \otimes b^q \in A^{\otimes q}$ by:

$$\begin{aligned}
\rho_{p,q}(\underline{a} \otimes \underline{b}) = &\ \varepsilon(a^1) a^2 \otimes \cdots \otimes a^p \otimes \underline{b} \\
&+ \sum_{i=1}^{p-1} (-1)^i a^1 \otimes \cdots \otimes a^i a^{i+1} \otimes \cdots \otimes a^p \otimes \underline{b} \\
&+ (-1)^p \sum_i \varepsilon(a_i^p) a^1 \otimes \cdots \otimes a^{p-1} \otimes \underline{b}_i \\
&+ (-1)^p \sum_j \varepsilon(b_j^1) \underline{a}_j \otimes b^2 \otimes \cdots \otimes b^q \\
&+ \sum_{i=1}^{q-1} (-1)^{i+p} \underline{a} \otimes b^1 \otimes \cdots \otimes b^i b^{i+1} \otimes \cdots \otimes b^q \\
&+ (-1)^{p+q} \varepsilon(b^q) \underline{a} \otimes b^1 \otimes \cdots \otimes b^{q-1},
\end{aligned}$$

where: $\sigma_q \ldots \sigma_1(a^p \otimes \underline{b}) = \sum_i \underline{b}_i \otimes a_i^p$, and for all $i$, $\underline{b}_i \in A^{\otimes q}$
$\sigma_1 \ldots \sigma_p(\underline{a} \otimes b^1) = \sum_j b_j^1 \otimes \underline{a}_j$, and for all $j$, $\underline{a}_j \in A^{\otimes p}$.

**Example:** Let $\sigma$ be the flip $\tau$, and $d_H$ the differential of the Hochschild complex of $A$ with coefficients in $\mathbb{K}$. Then a braided $\tau$-commutative algebra is just a commutative algebra, the quantum shuffle product is the shuffle product, and

$$\rho_{p,q} = d_H \otimes \operatorname{Id}^{\otimes q} + (-1)^p \operatorname{Id}^{\otimes p} \otimes d_H.$$

Thus, we find the classical result:

**Theorem 2.1** *Let $A$ be a commutative algebra. Then the Hochschild homology of $A$ with coefficients in $\mathbb{K}$ is a commutative graded algebra.*

As above, we can also consider normalized versions of the statements.

## 2.3 A new proof of the Hochschild-Serre identity.

Let $G$ be a group, $H$ a normal subgroup. If $N$ is a $G$-module we denote by $C^*(G, N)$ the chain complex of normalized cochains from $G$ to $N$. To determine a relation between the cohomology of $H$ and the one of $G$, Hochschild and Serre [7] have defined a map, later called Hochschild-Serre identity, between $C^*(G, N)$ and $C^*(G, C^*(H, N))$. They have shown, using combinatorial tricks, that this map is a morphism of complexes. This morphism was interpreted in terms of shuffles by N. Habegger, V. Jones, O. Pino Ortiz, J. Ratcliffe [6]. Using this framework, we extend this identity to Hopf algebras with invertible antipode. In particular, we show that the Hochschild-Serre identity is a dual statement of Corollary (2.1).

Let $A$ be a braided $\sigma$-commutative algebra and $M$ an $A$-bimodule. We have for all non negative integers $p, q$ a natural isomorphism:

$$\Lambda_{p,q} : \mathrm{Hom}_{\mathbb{K}}(A^{\otimes p} \otimes A^{\otimes q}, M) \simeq \mathrm{Hom}_{\mathbb{K}}(A^{\otimes q}, \mathrm{Hom}_{\mathbb{K}}(A^{\otimes p}, M)).$$

From this, we deduce the following identity of Hochschild-Serre type:

**Theorem 2.2** *There exists a differential $\rho$ on $\mathrm{Hom}_{\mathbb{K}}(T(A), Hom_{\mathbb{K}}(T(A), M))$ such that $\Lambda \circ {}^t\varphi_\sigma$ is a morphism of complexes from the Hochschild chain complex $\mathrm{Hom}_{\mathbb{K}}(T(A), M)$ to $(\mathrm{Hom}_{\mathbb{K}}(T(A), \mathrm{Hom}_{\mathbb{K}}(T(A), M)), \rho)$.*

**Proof:** The idea of the proof is the same as for the homological case. Instead of considering the tensor product $\widetilde{\otimes}$ we use the set, denoted by $\mathrm{Hom}_\sigma$, of the morphisms of $A \otimes_\sigma A$-bimodules.

In this framework the isomorphism $\Gamma_n$ (defined for all positive integer $n$) between $\mathrm{Hom}_\sigma(\beta_n(A \otimes_\sigma A, A \otimes_\sigma A), M)$ and $\mathrm{Hom}_{\mathbb{K}}((A \otimes_\sigma A)^{\otimes n}, M)$ arises naturally from the $A \otimes_\sigma A$-bimodule structure of $\beta_n(A \otimes_\sigma A, A \otimes_\sigma A)$.

For all positive integers $p, q$ the isomorphism $\Omega_{p,q}$ between $\mathrm{Hom}_\sigma(\beta_p(A, A) \otimes \beta_q(A, A), M)$ and $\mathrm{Hom}_{\mathbb{K}}(A^{\otimes p} \otimes A^{\otimes q}, M)$ is defined for $\underline{a} \in A^{\otimes p}, \underline{b} \in A^{\otimes q}$ by:

$$\Omega(f)(\underline{a} \otimes \underline{b}) = f(1_A[\underline{a}]1_A \otimes 1_A[\underline{b}]1_A).$$

The inverse is defined for $a[\underline{a}]a' \in \beta_p(A, A), b[\underline{b}]b' \in \beta_q(A, A)$ by:

$$\Omega^{-1}(f)(a'[\underline{a}]a'' \otimes b'[\underline{b}]b'') = \sum_i a' w_i \cdot f(\underline{w_i} \otimes \underline{\tau_i}) \cdot \tau_i b'',$$

with $\sigma_{p+q+1} \ldots \sigma_{p+2} \sigma_1 \ldots \sigma_{p+1}(\underline{a} \otimes a'' \otimes b' \otimes \underline{b}) = \sum_i w_i \otimes \underline{w_i} \otimes \underline{\tau_i} \otimes \tau_i$.

With these isomorphisms one can show that the map ${}^t\varphi_\sigma$ is a morphism of complexes. Using the isomorphism $\Lambda$ we obtain the result.

For all positive integers $p, q$, the differential $\rho$ is given on $f \in \mathrm{Hom}_{\mathbb{K}}(A^{\otimes p}, \mathrm{Hom}_{\mathbb{K}}(A^{\otimes q}, M))$ by:

**i)** for $a^1 \otimes \cdots \otimes a^{p+1} \in A^{\otimes p+1}$ and $\underline{b} \in A^{\otimes q}$:

$$\begin{aligned}
\rho(f)(a^1 \otimes \cdots \otimes a^{p+1})(\underline{b}) =\ & (-1)^q \sum_j a_j^1 . f(a^2 \otimes \cdots \otimes a^{p+1})(\underline{b}_j) \\
& + \sum_{i=1}^p (-1)^{i+q} f(a^1 \otimes \cdots \otimes a^i a^{i+1} \otimes \cdots \otimes a^{p+1})(\underline{b}) \\
& + (-1)^{p+q+1} f(a^1 \otimes \cdots \otimes a^p)(\underline{b}) . a^{p+1},
\end{aligned}$$

with $\sigma_1 \ldots \sigma_q(\underline{b} \otimes a^1) = \sum_j a_j^1 \otimes \underline{b}_j$, and for all $j$, $\underline{b}_j \in A^{\otimes q}$.

**ii)** for $\underline{a} \in A^{\otimes p}$ and $b^1 \otimes \cdots \otimes b^{q+1} \in A^{\otimes q+1}$:

$$\begin{aligned}
\rho(f)(\underline{a})(b^1 \otimes \cdots \otimes b^{q+1}) =\ & b^1 . f(\underline{a})(b^2 \otimes \cdots \otimes b^{q+1}) \\
& + \sum_{i=1}^q (-1)^i f(\underline{a})(b^1 \otimes \cdots \otimes b^i b^{i+1} \otimes \cdots \otimes b^{q+1}) \\
& + (-1)^{q+1} \sum_i f(\underline{a_i})(b^1 \otimes \cdots \otimes b^q) . b_i^{q+1},
\end{aligned}$$

with $\sigma_p \ldots \sigma_1(b^{q+1} \otimes \underline{a}) = \sum_i \underline{a}_i \otimes b_i^{q+1}$, and for all $i$, $\underline{a}_i \in A^{\otimes p}$.

□

As for the corollary 2.1, we show a normalized version of this result. For all non negative integer $p$, we denote by $\text{Hom}_N(A^{\otimes p}, M)$ the vector space of normalized maps of $\text{Hom}_\mathbb{K}(A^{\otimes p}, M)$, i.e. such that $f(a^1 \otimes \cdots \otimes a^p) = 0$ whenever one of the $a^i$ is equal to $1_A$. We have:

**Theorem 2.3** *There exists a differential $\rho$ on $\text{Hom}_N(T(A), Hom_N(T(A), M))$ such that $\Lambda \circ {}^t\varphi_\sigma$ induces a morphism of complexes between the normalized Hochschild chain complex $\text{Hom}_N(T(A), M)$ and $(\text{Hom}_N(T(A), \text{Hom}_N(T(A), M)), \rho)$.*

**Example:** In particular, this theorem can be applied to the braided $\sigma$-commutative algebra associated to a Hopf algebra with invertible antipode. Thus, we have a Hochschild-Serre identity for Hopf algebras with invertible antipode:

**Corollary 2.3** *Let $A$ be a Hopf algebra with invertible antipode and $H$ a sub-Hopf algebra stable for the adjoint action. The map $\Lambda \circ {}^t\varphi_\sigma$ induces a morphism of complexes between the normalized Hochschild chain complex and $(\text{Hom}_N(T(A), \text{Hom}_N(T(H), M)), \rho)$.*

In particular, for a group $G$, a normal subgroup $H$ and a left $\mathbb{K}[G]$-module $H$, the last corollary applied to $\mathbb{K}[G]$, $\mathbb{K}[H]$ and $M$ (which is a $\mathbb{K}[G]$-bimodule with the right action given by the counit of $\mathbb{K}[G]$) shows that the Hochschild-Serre identity is a statement dual to the result that the quantum shuffle product is a morphism of complexes.

## 3 The Hochschild-Serre identity in the multiplicative case.

In order to study the deformation by multiplicative coefficients of a class of Drinfeld Doubles, we prove a multiplicative version of the Hochschild-Serre identity. For this, we consider the abelian groups defined by M. Sweedler [17], $\text{Reg}(C, M)$, of invertible maps of $\text{Hom}_\mathbb{K}(C, M)$ for the convolution product, where $C$ is a cocommutative coalgebra and $M$ a commutative algebra. The convolution product is defined, for $f, g \in \text{Hom}_\mathbb{K}(C, M)$ by:

$$f * g = \mu_M(f \otimes g)\Delta_C.$$

We first give a multiplicative shuffle map for such abelian groups. As a consequence, we prove a multiplicative Hochschild-Serre identity for cocommutative Hopf algebras.

### 3.1 The multiplicative shuffle map.

In the additive case, we have considered $A \otimes_\sigma A$-bimodules. So we will work with invertible morphisms of $A \otimes_\sigma A$-bimodules. For this, we need complementary structures on $C$ and $M$.

**Definition 3.1** *Let $H$ be a bialgebra.*

*A coalgebra $C$, which is an $H$-bimodule, is called an $H$-bimodule coalgebra if $\Delta_C$ and $\varepsilon_C$ are $H$-bimodule morphisms.*

*An algebra $M$, which is an $H$-bimodule, is called an $H$-bimodule algebra if $\mu_M$ and $\eta_M$ are $H$-bimodule morphisms.*

If $C$ is an $H$-bimodule coalgebra, then for $a \in H$ and $c \in C$:

$$\Delta_C(a \cdot c) = \sum a_{(1)} \cdot c_{(1)} \otimes a_{(2)} \cdot c_{(2)}, \text{ and } \varepsilon_C(a \cdot c) = \varepsilon_H(a)\varepsilon_C(c).$$

If $M$ is an $H$-bimodule algebra, then for $a \in H$ and $\alpha, \beta \in M$:

$$a \cdot (\alpha\beta) = \sum (a_{(1)} \cdot \alpha)(a_{(2)} \cdot \beta), \text{ and } a \cdot 1_M = \varepsilon(a)1_M.$$

Thus, for $C$ a cocommutative $H$-bimodule coalgebra and $M$ a commutative $H$-bimodule algebra we can consider the subgroup $\text{Reg}_B(C, M)$ of the invertible $H$-bimodule morphisms from $C$ to $M$.

**Remark:** For $M$ fixed, $\text{Reg}_B(., M)$ is a contravariant functor between the category of cocommutative $H$-bimodule coalgebras and the category of abelian groups.

Let $A$ be a cocommutative Hopf algebra and $\sigma$ be the Woronowicz braiding. In Proposition 1.1 we have shown that $A \otimes_\sigma A$ is an algebra. The following lemma proves that it is a bialgebra.

**Lemma 3.1** *The Woronowicz braiding and its inverse are morphisms of coalgebras.*

Using the precedent lemma one can check that for all non negative integers $p, q$, the $A \otimes_\sigma A$-bimodule $\beta_p(A, A) \otimes \beta_q(A, A)$ is a cocommutative $A \otimes_\sigma A$-bimodule coalgebra, with the coproduct given by the tensor product of coalgebras.

Thus, if $M$ is a commutative $A \otimes_\sigma A$-bimodule algebra, we can consider the graded abelian groups:

$$\text{Reg}_B(\beta(A, A) \otimes \beta(A, A), M) = \bigoplus_{n \geq 0} \bigoplus_{p+q=n} \text{Reg}_B(\beta_p(A, A) \otimes \beta_q(A, A), M)$$

and

$$\text{Reg}_B(\beta(A, A) \times \beta(A, A), M) = \bigoplus_{n \geq 0} \text{Reg}_B(\beta_n(A, A) \otimes \beta_n(A, A), M).$$

In subsection 1.3, we have shown that the morphisms $\text{Id} \otimes s_i$, $s_i \otimes \text{Id}$, $\text{Id} \otimes d_i$ and $d_i \otimes \text{Id}$ are $A \otimes_\sigma A$-bimodule morphisms. Lemma 3.1 shows that they are coalgebra maps, so we can apply the functor $\text{Reg}_B(., M)$ to define:

$$\widetilde{s}_i^{[1]} = \text{Reg}_B(., M)(s_i \otimes \text{Id}), \qquad \widetilde{s}_i^{[2]} = \text{Reg}_B(., M)(\text{Id} \otimes s_i)$$
$$\widetilde{d}_i^{[1]} = \text{Reg}_B(., M)(d_i \otimes \text{Id}), \qquad \widetilde{d}_i^{[2]} = \text{Reg}_B(., M)(\text{Id} \otimes d_i)$$

The simplicial properties of $s_i$ and $d_i$ give analogue cosimplicial properties for $\widetilde{s}_i^{[l]}$, and $\widetilde{d}_i^{[l]}$, with $l \in \{1, 2\}$. So we have the dual multiplicative notion of tensor and cartesian product of simplicial chain complexes:

**Proposition 3.1** *The graded group* $\operatorname{Reg}_B(\beta(A,A) \otimes \beta(A,A), M)$ *is a chain complex. For all non negative integers* $p, q$, *the differential* $\widetilde{D_{p,q}}$ *from* $\operatorname{Reg}_B(\beta_p(A,A) \otimes \beta_q(A,A), M)$ *to* $\operatorname{Reg}_B(\beta_{p+1}(A,A) \otimes \beta_q(A,A), M) \oplus \operatorname{Reg}_B(\beta_p(A,A) \otimes \beta_{q+1}(A,A), M)$, *is defined by:*

$$\widetilde{D_{p,q}}(f) = \left( \overset{p+1}{\underset{i=0}{\ast}} \widetilde{d}_i^{[1]}(f^{(-1)^i}), \overset{q+1}{\underset{i=0}{\ast}} \widetilde{d}_i^{[2]}(f^{(-1)^{p+i}}) \right).$$

*The graded group* $\operatorname{Reg}_B(\beta(A,A) \times \beta(A,A), M)$ *is a chain complex with, for all non negative integer* $n$, *the differential:*

$$\widetilde{D_n} : \operatorname{Reg}_B(\beta_n(A,A) \otimes \beta_n(A,A), M) \to \operatorname{Reg}_B(\beta_{n+1}(A,A) \otimes \beta_{n+1}(A,A), M),$$

*defined by:*

$$\widetilde{D_n}(f) = \overset{n+1}{\underset{i=0}{\ast}} \widetilde{d}_i^{[1]} \circ \widetilde{d}_i^{[2]}(f^{(-1)^i}).$$

The multiplicative dual shuffle map is defined for all positive integers $p, q$ by:

$$\widetilde{g}_{p,q} : \operatorname{Reg}_B(\beta_{p+q}(A,A) \otimes \beta_{p+q}(A,A), M) \to \operatorname{Reg}_B(\beta_p(A,A) \otimes \beta_q(A,A), M),$$

and

$$\widetilde{g}_{p,q}(f) = \underset{\nu \in S_{p,p+q}}{\ast} \left( \widetilde{s}_{\nu(p+q)-1}^{[1]} \circ \cdots \circ \widetilde{s}_{\nu(p+1)-1}^{[1]} \right) \circ \left( \widetilde{s}_{\nu(p)-1}^{[2]} \circ \cdots \circ \widetilde{s}_{\nu(1)-1}^{[2]} \right)(f^{(-1)^{|\nu|}}).$$

In the additive case, the proof of Theorem 1.1 arises from the properties of the simplicial objects and of the $(p,q)$-shuffle. In the multiplicative dual case the $(p,q)$-shuffle properties hold, and we have the analogue properties for the simplicial objects, thus:

**Theorem 3.1** *The multiplicative dual shuffle map is a morphism of complexes.*

### 3.2 The multiplicative Hochschild-Serre identity.

Let $M$ be a commutative $A$-bimodule algebra. It is clear that $M$ is a commutative $A \otimes_\sigma A$-bimodule algebra.

Using Lemma 3.1, we can apply the functor $\operatorname{Reg}_B(.,M)$ to the maps $\sigma_i$. So for all non negative integer $n$ we define the dual of $\Psi_n$, and for $p \geq 2$ and all $w \in \Sigma_p$ the dual of $T_\sigma(w)$:

$$\widetilde{\Psi}_n : \operatorname{Reg}_B(\beta_n(A \otimes_\sigma A, A \otimes_\sigma A), M) \to \operatorname{Reg}_B(\beta_n(A,A) \otimes \beta_n(A,A), M),$$
$$\widetilde{T_\sigma(w)} : \operatorname{Reg}(A^{\otimes p}, M) \to \operatorname{Reg}(A^{\otimes p}, M).$$

As in the additive case, $\widetilde{\Psi}$ is a morphism of complexes. So we have a morphism of complexes from $\operatorname{Reg}_B(\beta(A \otimes_\sigma A, A \otimes_\sigma A), M)$ to $\operatorname{Reg}_B(\beta(A,A) \otimes \beta(A,A), M)$.

Furthermore, the isomorphisms of section 2.3, can be written multiplicatively i.e. for all positive integers $p, q$:

$$\operatorname{Reg}_B(\beta_p(A \otimes_\sigma A, A \otimes_\sigma A), M) \simeq \operatorname{Reg}\left((A \otimes_\sigma A)^{\otimes p}, M\right),$$
$$\operatorname{Reg}_B(\beta_p(A,A) \otimes \beta_q(A,A), M) \simeq \operatorname{Reg}(A^{\otimes p} \otimes A^{\otimes q}, M),$$
$$\operatorname{Reg}(A^{\otimes p} \otimes A^{\otimes q}, M) \overset{\widetilde{\Lambda}}{\simeq} \operatorname{Reg}\left(A^{\otimes q}, \operatorname{Hom}_{\mathbb{K}}(A^{\otimes p}, M)\right).$$

As in the additive case, one can prove:

**Theorem 3.2** *The multiplicative dual quantum shuffle product map:*

$$\widetilde{\varphi_\sigma} : \bigoplus_{n \geq 0} \operatorname{Reg}(A^{\otimes n}, M) \to \bigoplus_{n \geq 0} \bigoplus_{p+q=n} \operatorname{Reg}(A^{\otimes p} \otimes A^{\otimes q}, M),$$

*defined for $f \in \operatorname{Reg}(A^{\otimes n}, M)$, $\underline{a} \in A^{\otimes p}$ and $\underline{b} \in A^{\otimes n-p}$ by:*

$$\widetilde{\varphi_\sigma}(f)(\underline{a} \otimes \underline{b}) = \underset{w \in S_{p,n}}{\ast} \widetilde{T_\sigma(w)}(f^{(-1)^{|w|}})(\underline{a} \otimes \underline{b}),$$

*is a morphism of complexes.*

Using the isomorphism $\widetilde{\Lambda}$, we have a morphism of complexes:

$$\widetilde{\Lambda} \circ \widetilde{\varphi_\sigma} : \left( \bigoplus_{n \geq 0} \operatorname{Reg}(A^{\otimes n}, M), \widetilde{\delta} \right) \longrightarrow \left( \bigoplus_{n \geq 0} \bigoplus_{p+q=n} \operatorname{Reg}(A^{\otimes q}, \operatorname{Hom}_{\mathbb{K}}(A^{\otimes p}, M)), \widetilde{\rho} \right).$$

The first chain complex is the multiplicative Hochschild chain complex. The second chain complex can be understood as the total chain complex of the following chain bicomplex:

for all positive integers $p, q$, put $C^{p,q}(A, M) = \operatorname{Reg}(A^{\otimes p}, \operatorname{Hom}_{\mathbb{K}}(A^{\otimes q}, M))$,

the vertical differential is given by the map $\widetilde{D}_p^{(q)} : C^{p,q}(A, M) \to C^{p+1,q}(A, M)$, defined for $f \in C^{p,q}(A, M)$, $\underline{a} = a^1 \otimes \cdots \otimes a^{p+1} \in A^{\otimes p+1}$ and $\underline{x} \in A^{\otimes q}$ by:

$$\begin{aligned}
\widetilde{D}_p^{(q)}(f)(\underline{a})(\underline{x}) = & \sum a^1_{(1)} \cdot f^{(-1)^q}(a^2_{(1)} \otimes \cdots \otimes a^{p+1}_{(1)})(\mathrm{ad}_{a^1_{(2)}} \underline{x}_{(1)}) \\
& \prod_{i=1}^{p} f^{(-1)^{q+i}}(a^1_{(i+2)} \otimes \cdots \otimes a^i_{(i+1)} a^{i+1}_{(i+1)} \otimes \cdots \otimes a^{p+1}_{(i+1)})(\underline{x}_{(i+1)}) \\
& f^{(-1)^{q+p+1}}(a^1_{(p+3)} \otimes \cdots \otimes a^p_{(p+2)})(\underline{x}_{(p+2)}) \cdot a^{p+1}_{(p+2)}
\end{aligned}$$

where $\mathrm{ad}_{a^1} \underline{x} = \sum \mathrm{ad}_{a^1_{(1)}} x^1 \otimes \cdots \otimes \mathrm{ad}_{a^1_{(q)}} x^q$.

the horizontal differential is given by the map $\widetilde{d}_q^{(p)} : C^{p,q}(A, M) \to C^{p,q+1}(A, M)$, defined for $f \in C^{p,q}(A, M)$, $\underline{a} \in A^{\otimes p}$ and $\underline{x} = x^1 \otimes \cdots \otimes x^{q+1} \in A^{\otimes q+1}$ by:

$$\begin{aligned}
\widetilde{d}_q^{(p)}(f)(\underline{a})(\underline{x}) = & \sum x^1_{(1)} \cdot f(\underline{a}_{(1)})(x^2_{(1)} \otimes \cdots \otimes x^{q+1}_{(1)}) \\
& \prod_{i=1}^{q} f^{(-1)^i}(\underline{a}_{(i+1)})(x^1_{(i+1)} \otimes \cdots \otimes x^i_{(i+1)} x^{i+1}_{(i+1)} \otimes \cdots \otimes x^{q+1}_{(i+1)}) \\
& f^{(-1)^{q+1}}(\underline{a}_{(q+2)})(x^1_{(q+2)} \otimes \cdots \otimes x^q_{(q+2)}) \cdot \mathrm{ad}_{\underline{a}_{(q+3)}} x^{q+1}
\end{aligned}$$

where $\mathrm{ad}_{\underline{a}} x^{q+1} = \mathrm{ad}_{a^1 \ldots a^p} x^{q+1}$.

With this definition the differential $\widetilde{\rho}$ is given, for $f \in C^{p,q}(A, M)$, by:

$$\widetilde{\rho}(f) = (\widetilde{D}_p^{(q)}(f), \widetilde{d}_q^{(p)}(f)).$$

From the construction of the differential $\widetilde{\rho}$, for all non negative integer $n$ we have that $\left( C^{\cdot,n}(A, M), \widetilde{D}_\cdot^{(n)} \right)$ and $\left( C^{n,\cdot}(A, M), \widetilde{d}_\cdot^{(n)} \right)$ are chain complexes. Moreover, the differentials are compatible in the chain bicomplex sense, i.e. for $f \in C^{p,q}(A, M)$:

$$\left( \widetilde{D}_p^{(q+1)} \circ \widetilde{d}_q^{(p)} \right)(f) \ast \left( \widetilde{d}_q^{(p+1)} \circ \widetilde{D}_p^{(q)} \right)(f) = \eta_M \circ \varepsilon_A.$$

**Remark:** As in the additive case, we have a normalized version of this result.

For all non negative integer $n$, denote by $\mathrm{Hom}_N(A^{\otimes n}, M)$ the set of all normalized maps, i.e. the maps $f \in \mathrm{Hom}_{\mathbb{K}}(A^{\otimes n}, M)$ such that $f(a^1 \otimes \cdots \otimes a^n) = \varepsilon(a^1 \otimes \cdots \otimes a^n)$ whenever one of the $a^i$ is equal to $1_A$.

For all non negative integers $p, q$, put $C_N^{p,q}(A, M) = \mathrm{Reg}_N(A^{\otimes p}, \mathrm{Hom}_N(A^{\otimes q}, M))$. With this notation, we have:

**Corollary 3.1** *Let $A$ be a cocommutative Hopf algebra.*
*The map:*

$$\widetilde{\Lambda} \circ \widetilde{\varphi_\sigma} : \left( \bigoplus_{n \geq 0} \mathrm{Reg}_N(A^{\otimes n}, M), \widetilde{\delta} \right) \longrightarrow \left( \bigoplus_{n \geq 0} \bigoplus_{p+q=n} C_N^{p,q}(A, M), \widetilde{\rho} \right),$$

*is a morphism of complexes.*

### 3.3 Example in degree 3.

Let $f \in \mathrm{Reg}_N(A^{\otimes 3}, M)$ and $x \otimes y \otimes z \otimes t \in A^{\otimes 4}$. Then:

$$\widetilde{\delta}_3(f)(x \otimes y \otimes z \otimes t) = \sum (x_{(1)} \cdot f(y_{(1)} \otimes z_{(1)} \otimes t_{(1)})) f^{-1}(x_{(2)} y_{(2)} \otimes z_{(2)} \otimes t_{(2)})$$
$$f(x_{(3)} \otimes y_{(3)} z_{(3)} \otimes t_{(3)}) f^{-1}(x_{(4)} \otimes y_{(4)} \otimes z_{(4)} t_{(4)})$$
$$(f(x_{(5)} \otimes y_{(5)} \otimes z_{(5)}) \cdot t_{(5)})$$

The map $\widetilde{\Lambda} \circ \widetilde{\varphi_\sigma}$ sends $f$ to $(f_0, f_1, f_2, f_3)$, where for $l \in \{1, 2, 3\}$:

$$f_l \in \mathrm{Reg}_N(A^{\otimes l}, \mathrm{Hom}_N(A^{\otimes 3-l}, M)).$$

We have $f_0 = f_3 = f$, and for all elements $x, y, a, b$ of $A$:

$$f_1(x)(a \otimes b) = \sum f(a_{(1)} \otimes b_{(1)} \otimes x_{(1)}) f^{-1}(a_{(2)} \otimes x_{(2)} \otimes \mathrm{ad}_{x_{(3)}} b_{(2)}) f(x_{(4)} \otimes \mathrm{ad}_{x_{(5)}} a_{(3)} \otimes \mathrm{ad}_{x_{(6)}} b_{(3)}),$$

$$f_2(x \otimes y)(a) = \sum f(a_{(1)} \otimes x_{(1)} \otimes y_{(1)}) f^{-1}(x_{(2)} \otimes \mathrm{ad}_{x_{(3)}} a_{(2)} \otimes y_{(2)}) f(x_{(4)} \otimes y_{(3)} \otimes \mathrm{ad}_{x_{(5)} y_{(4)}} a_{(3)}).$$

Let $(f_0, f_1, f_2, f_3) \in \bigoplus_{p+q=3} C_N^{p,q}(A, M)$. Then:

$$\widetilde{\rho}(f_0, f_1, f_2, f_3) = \left( \widetilde{\delta}_3(f_0), \widetilde{D}_0^{(3)}(f_0) * \widetilde{d}_2^{(1)}(f_1), \widetilde{D}_1^{(2)}(f_1) * \widetilde{d}_1^{(2)}(f_2), \widetilde{D}_2^{(1)}(f_2) * \widetilde{d}_0^{(3)}(f_3), \widetilde{\delta}_3(f_3) \right).$$

## 4 Quasi-Hopf algebras obtained from Drinfeld Doubles.

In order to explain a construction of quasi-Hopf algebras, associated to finite dimensional cocommutative Hopf algebras, we will use the multiplicative Hochschild-Serre identity. This type of construction was first considered by R. Dijkgraaf, V. Pasquier and P. Roche [3] for the Hopf algebra of a finite group and by D. Bulacu and F. Panaite [2] for finite dimensional cocommutative Hopf algebras.

They deform the product and the coproduct of the Drinfeld double by two particular maps obtained from a normalized 3-cocycle of the multiplicative Hochschild chain complex. Then they prove that the deformed Drinfeld double is a quasi-Hopf algebra. The choice of these two maps was only motivated by the fact that they verify the conditions allowing to build a quasi-Hopf algebra.

We will first extend the construction to a more general deformation and explain the link between the deformation and the normalized 3-cocycle of the total complex. At last, using the multiplicative Hochschild-Serre identity, we clarify the choice of R. Dijkgraaf, V. Pasquier, P. Roche and D. Bulacu, F. Panaite.

## 4.1 Definition of the quasi-Hopf algebras.

In the representation theory of Hopf algebras, it is well known that the coproduct induces a tensor product of representations which is strictly associative. This is a consequence of the coproduct coassociativity.

A way to have a theory with a non associative tensor product of representations is to use representations of a quasi-Hopf algebra. These objects were first defined by V.G. Drinfeld [4]:

**Definition 4.1** *An algebra $A$ is a quasi-Hopf algebra if there exist morphisms of algebras $\Delta : A \to A \otimes A$ and $\varepsilon : A \to \mathbb{K}$, an anti-automorphism of algebras $S : A \to A$, an invertible element $\phi \in A \otimes A \otimes A$ and elements $\alpha, \beta \in A$ such that:*

*i) $\Delta$ is quasi-coassociative, i.e. for all $a \in A$:*

$$(\mathrm{Id} \otimes \Delta)\Delta(a) = \phi(\Delta \otimes \mathrm{Id})\Delta(a)\phi^{-1},$$

*ii) the associator $\phi$ satisfies the pentagonal equality:*

$$(\mathrm{Id} \otimes \mathrm{Id} \otimes \Delta)(\phi) \cdot (\Delta \otimes \mathrm{Id} \otimes \mathrm{Id})(\phi) = (1_A \otimes \phi) \cdot (\mathrm{Id} \otimes \Delta \otimes \mathrm{Id})(\phi) \cdot (\phi \otimes 1_A),$$

*iii) the counit $\varepsilon$ satisfies:*

$$(\varepsilon \otimes \mathrm{Id})\Delta = \mathrm{Id} = (\mathrm{Id} \otimes \varepsilon)\Delta, \qquad (\mathrm{Id} \otimes \varepsilon \otimes \mathrm{Id})(\phi) = 1_A \otimes 1_A, \qquad (5)$$

*iv) for all $a \in A$:*

$$\sum S(a_{(1)})\alpha a_{(2)} = \varepsilon(a)\alpha, \qquad \sum a_{(1)}\beta S(a_{(2)}) = \varepsilon(a)\beta, \qquad (6)$$

*v) if $\phi = \sum_i X_i \otimes Y_i \otimes Z_i$ and $\phi^{-1} = \sum_i P_i \otimes Q_i \otimes R_i$ then:*

$$\sum_i X_i \beta S(Y_i)\alpha Z_i = 1_A, \qquad \sum_i S(P_i)\alpha Q_i \beta S(R_i) = 1_A. \qquad (7)$$

## 4.2 Deformation of the Drinfeld Double into a Hopf algebra.

Let $A$ be a finite dimensional cocommutative Hopf algebra.

The space $A^*$ is an $A$-bimodule with the actions defined for $f \in A^*$ and $x, a \in A$ by:

$$(xf)(a) = f(ax), \qquad (fx)(a) = f(xa).$$

The Drinfeld double of $A$, denoted by $\mathcal{D}(A)$, is the Hopf algebra defined as a vector space by $A^* \otimes A$. For $f, g \in A^*$ and $x, y \in A$, the structure is given by:

$$\begin{aligned}
(f \otimes x) \cdot (g \otimes y) &= \sum f(x_{(1)} g \mathrm{S}(x_{(2)})) \otimes x_{(3)} y, \\
\Delta(f \otimes x) &= \sum (f_{(1)} \otimes x_{(1)}) \otimes (f_{(2)} \otimes x_{(2)}), \\
\mathrm{S}(g \otimes x) &= \sum \mathrm{S}(x_{(1)}) \mathrm{S}(g) x_{(2)} \otimes \mathrm{S}(x_{(3)}).
\end{aligned}$$

The unit and counit are those given by the tensor product of bialgebras.

Let $\theta \in C_N^{2,1}(A, \mathbb{K})$, and $\gamma \in C_N^{1,2}(A, \mathbb{K})$. As $A$ is a finite dimensional vector space, we have $(A \otimes A)^* \simeq A^* \otimes A^*$. So for $x \in A$, we have $\gamma(x) \in A^* \otimes A^*$, and we denote by $\gamma(x) = \sum \gamma(x)_{(1)} \otimes \gamma(x)_{(2)}$.

The product and the coproduct are deformed by:

$$(f \otimes x) \cdot (g \otimes y) = \sum f(x_{(1)} g \mathrm{S}(x_{(2)})) \theta(x_{(3)} \otimes y_{(1)}) \otimes x_{(4)} y_{(2)}, \tag{8}$$

$$\Delta_\gamma(f \otimes x) = \sum (f_{(1)} \gamma(x_{(1)})_{(1)} \otimes x_{(2)}) \otimes (f_{(2)} \gamma(x_{(1)})_{(2)} \otimes x_{(3)}). \tag{9}$$

**Remark:** R. Dijkgraaf, V. Pasquier and P. Roche [3], and D. Bulacu and F. Panaite [2] have studied such deformations associated to a particular choice of $\theta$ and $\gamma$.

More generally, I. Hofstetter [8] has studied Hopf algebras obtained by tensor product of an abelian matched pair of Hopf algebras. Our deformation is a particular case of this work in the following sense: the precedent braided structures are those obtained from the Hofstetter theory applied to the abelian matched pair $(A, A^*)$.

Thus, we have:

**Theorem 4.1 (Hofstetter)** *Let $A$ be a finite dimensional cocommutative Hopf algebra. The space $\mathcal{D}(A)$ is a Hopf algebra with the crossed product and coproduct if and only if $(\varepsilon, \gamma, \theta, \varepsilon)$ is a normalized $3$-cocycle of the total chain complex associated to $(C^{\cdot,\cdot}(A, \mathbb{K}), \widetilde{\rho})$. The unit and the counit are those given by the tensor product of bialgebras. The antipode is given, for $f \in A^*$ and $x \in A$, by:*

$$\mathrm{S}(f \otimes x) = \sum \left( \varepsilon \otimes \mathrm{S}(x_{(1)}) \right) \cdot \left( \theta^{-1}(x_{(2)} \otimes \mathrm{S}(x_{(3)})) \mathrm{S} \left( f \gamma^{-1}(x_{(4)})_{(1)} \right) \gamma^{-1}(x_{(4)})_{(2)} \otimes 1_A \right).$$

*We denote this Hopf algebra by $\mathcal{D}^{\theta,\gamma}(A)$.*

**Proof:** In her article Hofstetter proves the equivalent condition: $(\gamma, \theta)$ is a 2-cocycle of the subcomplex $\left( \bigoplus_{n \geq 1} \bigoplus_{p+q=n} C_N^{p,q}(A, \mathbb{K}), \widetilde{\rho} \right)$.

In our case, the associativity is equivalent to $\widetilde{D}_2^{(1)}(\theta) = \varepsilon$.

The coassociativity is equivalent to $\widetilde{d}_2^{(1)}(\gamma) = \varepsilon$.

The coproduct $\Delta_\gamma$ is a morphism of algebras if and only if $\widetilde{D}_1^{(2)}(\gamma) * \widetilde{d}_1^{(2)}(\theta) = \varepsilon$. □

## 4.3 Deformation of the Drinfeld Double into a quasi-Hopf algebra.

Instead of considering normalized 3-cocycle of the type $(\varepsilon, \gamma, \theta, \varepsilon)$, we consider normalized 3-cocycle of the type $(w, \gamma, \theta, \varepsilon)$. We will show that $\mathcal{D}^{\theta,\gamma}(A)$ is a quasi-Hopf algebra and that $w$ corresponds to the associator.

Let $(w, \gamma, \theta, \varepsilon)$ be a normalized 3-cocycle.

Using the properties of the cocycle and Theorem 4.1, it is easy to see that with the braided product (8) the space $\mathcal{D}^{\theta,\gamma}(A)$ is an algebra, and that the braided coproduct (9) and the counit are morphisms of algebras.

But the coproduct is not coassociative. This non-coassociativity is described by $w$.

In fact, the map $w$ is an element of $\text{Reg}_N(A^{\otimes 3}, \mathbb{K})$, so it can be viewed as an element of $(A^*)^{\otimes 3}$. We denote it by $w = \sum w_{(1)} \otimes w_{(2)} \otimes w_{(3)}$, and its inverse by $w^{-1} = \sum w_{(1)}^{-1} \otimes w_{(2)}^{-1} \otimes w_{(3)}^{-1}$.

Thus, we can define an invertible element $\phi \in \mathcal{D}^{\theta,\gamma}(A)^{\otimes 3}$ by:

$$\phi = \sum (w_{(1)}^{-1} \otimes 1_A) \otimes (w_{(2)}^{-1} \otimes 1_A) \otimes (w_{(3)}^{-1} \otimes 1_A).$$

The inverse is given by $\phi^{-1} = \sum (w_{(1)} \otimes 1_A) \otimes (w_{(2)} \otimes 1_A) \otimes (w_{(3)} \otimes 1_A)$.

With these notations, the cocycle condition $\widetilde{D}_0^{(3)}(w) * \widetilde{d}_2^{(1)}(\gamma) = \varepsilon$ is equivalent, for all $f \in A^*$ and $x \in A$, to:

$$(\text{Id} \otimes \Delta_\gamma)\Delta_\gamma(f \otimes x) = \phi (\Delta_\gamma \otimes \text{Id})\Delta_\gamma(f \otimes x)\phi^{-1}.$$

The cocycle condition $\widetilde{d}_3^{(0)}(w) = \varepsilon$ is equivalent to:

$$(\text{Id} \otimes \text{Id} \otimes \Delta_\gamma)(\phi) \cdot (\Delta_\gamma \otimes \text{Id} \otimes \text{Id})(\phi) = ((\varepsilon \otimes 1_A) \otimes \phi) \cdot (\text{Id} \otimes \Delta_\gamma \otimes \text{Id})(\phi) \cdot (\phi \otimes (\varepsilon \otimes 1_A)).$$

Thus, as suggested by remark 3.2 of [2]:

**Theorem 4.2** *Let $A$ be a finite dimensional cocommutative Hopf algebra.*
*Let $(w, \gamma, \theta, \varepsilon) \in \bigoplus_{p+q=3} C_N^{p,q}(A, \mathbb{K})$, $\phi = \sum(w_{(1)}^{-1} \otimes 1_A) \otimes (w_{(2)}^{-1} \otimes 1_A) \otimes (w_{(3)}^{-1} \otimes 1_A)$, $\alpha = \varepsilon \otimes 1_A$ and $\beta = \sum(w_{(1)} S(w_{(2)}) w_{(3)}) \otimes 1_A$.*
*The space $\mathcal{D}(A)$ is a quasi-Hopf algebra (denoted by $\mathcal{D}^{w,\theta,\gamma}(A)$) with the braided product (8), the coproduct (9), the unit, counit, and antipode of Theorem 4.1, the associator $\phi$ and the elements $\alpha, \beta$ if and only if $(w, \gamma, \theta, \varepsilon)$ is a normalized 3-cocycle of the total chain complex associated to $(C^{\cdots}(A, \mathbb{K}), \widetilde{\rho})$.*

**Proof:** Using Theorem 4.1 and the above discussion, we only have to prove the equalities (5), (6) and (7) of quasi-Hopf algebras.

Equations (5) arise from to the normalization of the cocycle.

Equations (6) follow from:

$$\widetilde{D}_2^{(1)}(\theta)(\sum S(x_{(1)}) \otimes x_{(2)} \otimes S(x_{(3)})) = \varepsilon(x)\varepsilon,$$
$$\widetilde{d}_2^{(1)}(\gamma) * D_0^{(3)}(w)(x)(\sum S(z_{(1)}) \otimes z_{(2)} \otimes S(z_{(3)}))) = \varepsilon(x)\varepsilon(z),$$

where $z, x \in A$.

Equations (7) are a consequence of the definitions of $\phi$, $\alpha$, $\beta$ and of the equality, given for $z \in A$ by:

$$\widetilde{\delta}_3(w)\left(\sum S(z_{(1)}) \otimes z_{(2)} \otimes S(z_{(3)}) \otimes z_{(4)}\right) = \varepsilon(z).$$

□

## 4.4 A class of examples explained by the Hochschild-Serre identity.

In the papers [3] of R. Dijkgraaf, V. Pasquier, P. Roche and [2] of D. Bulacu, F. Panaite, the authors make the precedent construction for a normalized 3-cocycle $w$ of the multiplicative Hochschild chain complex and the following particular maps, defined for $x, y, z \in A$ by:

$$\theta(x \otimes y)(z) = \sum w(z_{(1)} \otimes x_{(1)} \otimes \otimes y_{(1)})w^{-1}(x_{(2)} \otimes \mathrm{ad}_{x_{(3)}} z_{(2)} \otimes y_{(2)})$$
$$w(x_{(4)} \otimes y_{(3)} \otimes \mathrm{ad}_{x_{(5)}y_{(4)}} z_{(3)}),$$
$$\gamma(z)(x \otimes y) = \sum w(x_{(1)} \otimes y_{(1)} \otimes z_{(1)})w^{-1}(x_{(2)} \otimes z_{(2)} \otimes \mathrm{ad}_{z_{(3)}} y_{(2)})$$
$$w(z_{(4)} \otimes \mathrm{ad}_{z_{(5)}} x_{(3)} \otimes \mathrm{ad}_{z_{(6)}} y_{(3)}).$$

This non-trivial choice is explained by:

**Theorem 4.3** *Let $w$ be a normalized 3-cocycle of the multiplicative Hochschild chain complex. The constructions of R. Dijkgraaf, V. Pasquier, P. Roche and D. Bulacu, F. Panaite are those obtained by deforming the Drinfeld double with the image of $w$ under the multiplicative Hochschild-Serre identity.*

**Proof:** First, remark that using the results of section 3.3 the maps $\gamma$ and $\theta$ are exactly the maps, of $C_N^{1,2}(A, \mathbb{K})$ and $C_N^{2,1}(A, \mathbb{K})$, obtained by the Hochschild-Serre identity applied to $w$:

$$\widetilde{\Lambda} \circ \widetilde{\varphi_\sigma}(w) = (w, \gamma, \theta, w).$$

By Corollary 3.1, $(w, \gamma, \theta, w)$ is a normalized 3-cocycle of the total complex associated to $(C^{\cdots}(A, \mathbb{K}), \widetilde{\rho})$.

Furthermore:

**Lemma 4.1** *Let $w_1$ be a normalized 3-cocycle for the multiplicative Hochschild chain complex. Then:*

$(w, \gamma, \theta, \varepsilon)$ *is a normalized* 3-*cocycle* $\iff$ $(w, \gamma, \theta, w_1)$ *is a normalized* 3-*cocycle.*

**Proof:** Let $f \in \mathrm{Reg}_N(A^{\otimes p}, \mathbb{K})$. Then $\widetilde{d}_0^{(p)}(f) = \varepsilon$ and $\widetilde{D}_p^{(0)}(f) = \widetilde{d}_p^{(0)}(f) = \widetilde{\delta}_p(f)$. From this, we deduce immediately the result. $\square$

In particular, for $w_1 = w$, we can apply Theorem 4.2 to the normalized 3-cocycle $(w, \gamma, \theta, w)$. $\square$